\documentclass[12pt]{article}
\usepackage[utf8]{inputenc}
\usepackage{xcolor,soul}
\usepackage{amsmath,amssymb,amsthm,enumerate,graphicx}
\usepackage{mathtools}
\usepackage{ytableau}
\usepackage{tikz}
\tikzstyle{vertex}=[circle, draw, fill=black, inner sep=0pt, minimum size=4pt]
\usepackage{soul}

\usepackage{verbatim}
\usepackage{overpic}
\usepackage{geometry}

\usepackage{graphicx}

\theoremstyle{plain}
\newtheorem{theorem}{Theorem}
\newtheorem{proposition}[theorem]{Proposition}

\newtheorem{conjecture}[theorem]{Conjecture}
\theoremstyle{definition}

\newtheorem{definition}[theorem]{Definition}

\newcommand{\match}{\mathcal{M}}

\title{Matching Complexes of Outerplanar Graphs}
\author{Margaret Bayer\\Department of Mathematics\\University of Kansas\\Lawrence, Kansas, U.S.A.\\bayer@ku.edu \and Marija Jeli\'{c} Milutinovi\'{c}\\Faculty of Mathematics\\University of Belgrade\\ Belgrade, Serbia\\marijaj@matf.bg.ac.rs \and  Julianne Vega\thanks{This article is based on work supported by the National Science Foundation under Grant No.\ DMS-1440140 while the authors participated in the 2020/2021 Summer Research in Mathematics Program of the Mathematical Sciences Research Institute, Berkeley, California. M.\ Jeli\'{c} Milutinovi\'{c} has been supported by the Project No.\ 7744592 MEGIC
”Integrability and Extremal Problems in Mechanics, Geometry and
Combinatorics” of the Science Fund of Serbia, and by the Faculty of Mathematics University of Belgrade through the grant (No.\ 451-03-68/2022-14/200104) by the Ministry of
Education, Science, and Technological Development of the Republic of Serbia.}
\\Department of Mathematics\\Maret School\\Washington, D.C. U.S.A.\\jvega@maret.org}

\begin{document}
\maketitle

\begin{abstract}
An outerplanar graph is a planar graph that has a planar drawing with all vertices on the unbounded face.  The matching complex of a graph is the simplicial complex whose faces are subsets of disjoint edges of the graph.  In this paper we prove that the matching complexes of outerplanar graphs are contractible or homotopy equivalent to a wedge of spheres.  This extends known results about trees and polygonal line tilings.
\end{abstract}

\section{Introduction}

For a graph $G$, a matching is a set of pairwise disjoint edges, and the matching complex, $\match(G)$, is the simplicial complex whose faces are matchings of $G$.  (For any graph theory terminology, the reader may refer to \cite{West}.)
There has been much interest in the topology of matching complexes.  Most relevant to this paper are papers by Matsushita \cite{takahiro}, Jeli\'c Milutinovi\'c, et al.\ \cite{JelicEtAl}, and the previous papers of the current authors \cite{BJV,bjv-perfect}.  In addition, a recent paper of Gupta, et al.\ \cite{gupta2024matching} studies the matching complexes of categorical products of paths.
Also of note are the studies of Gorenstein, Cohen--Macaulay and Buchsbaum properties of matching complexes, by Nikseresht \cite{nikseresht} and Goeckner, et al.\ \cite{Goeckner-Buchsbaum}.  For a survey of earlier work, see Wachs \cite{wachs}.

The tools to study these complexes vary. Donovan and Scoville
\cite{donovan-scoville}
use star clusters and collapses (following Barmak \cite{barmak}) to provide alternative proofs for the homotopy type of the matching complexes for paths and cycles as well as to characterize the homotopy type of Dutch Windmill graphs as wedges of spheres.
There has also been progress using discrete Morse theory to study the topology of matching complexes for complete graphs; see Mondal, Mukherjee and Saha \cite{mondal2024topology}.
In the current paper, we will use tools for studying independence complexes developed by Adamaszek and Engstr\"om to characterize the matching complex of outerplanar graphs as wedges of spheres.

\section{Main Result}
\begin{definition}
A graph is {\em outerplanar} if and only if it has an embedding in the
plane so that all of its vertices lie on the unbounded face.
\end{definition}
A trivial example of an outerplanar graph is an {\em empty graph}, a graph with at least one vertex, but no edges.  For such a graph, the matching complex is
empty; in the following we assume all our graphs are nonempty.
\begin{theorem}\label{main-th}
If $G$ is a (nonempty) outerplanar graph, then the matching complex $\match(G)$ is
contractible or homotopy equivalent to a wedge of spheres.
\end{theorem}

The proof relies on several theorems that were stated and proved for independence complexes. For a graph $G$, the independence complex of $G$ is the simplicial complex whose faces are independent sets of vertices of $G$.  We translate these into results for matching complexes, using the following.  
\begin{definition}
The {\em line graph} $L(G)$ of a graph $G$ is the graph whose vertex set is
the set of edges of $G$ and whose edge set is the set of pairs of edges of $G$
that share a vertex.
\end{definition}

The following statement follows directly from the definitions.
\begin{proposition} \label{prop:math_ind}
The matching complex of $G$ is the independence complex of $L(G)$.
\end{proposition}

Define the (open) edge neighborhood $\operatorname{EN}_G(e)$
of an edge $e$ in the graph $G$ to be the
set of edges adjacent to $e$, and
the closed edge neighborhood of $e$ to be $\operatorname{EN}_G[e] = \operatorname{EN}_G(e) \cup \{e\}$. When the graph $G$ is clear from context we write $\operatorname{EN}(e)$ and $\operatorname{EN}[e]$, respectively.

We use the following (translated) theorems. 

\begin{proposition}[Adamaszek \cite{adamaszek-split}, Theorem~3.3]\label{prop:closed_neigh}
Let $G$ be a graph that contains two different edges $e$ and $h$ such that $\operatorname{EN}[e] \subseteq \operatorname{EN}[h]$. Then
\[\match(G) \simeq \match(G  \setminus\{h \}) \vee \Sigma\match(G  \setminus \operatorname{EN}[h]).\]
\end{proposition}

Notice if $G$ has at least two edges, and the closed neighborhood $\operatorname{EN}[h]$ consists of all edges of $G$,  then for any edge $e$, $\operatorname{EN}[e]\subseteq \operatorname{EN}[h]$. 
So $G\setminus \operatorname{EN}[h]$ is a graph with a nonempty vertex set but no edges.
The matching complex, which is a complex on the set of edges, thus has no
faces; we call it the void complex.
But the matching complex of $G$ includes a 0-dimensional face
corresponding to the edge $h$.  For the above proposition to apply to this
case, we adopt the convention that the suspension of the void complex is the 0-sphere $S^0$. 

\begin{proposition}[Engstr\"om \cite{engstrom2}, Lemma 2.4]\label{prop:open_neighborhood}
Let $G$ be a graph that contains two different edges $e$ and $h$ such that $\operatorname{EN}(e) \subseteq \operatorname{EN}(h)$. Then $\match(G)$ collapses to  $\match(G \setminus\{h\})$. Thus, $\match(G)\simeq \match(G \setminus\{h\})$.
\end{proposition}

\begin{proposition}[Engstr\"om \cite{engstrom-KM}, Lemma~2.2]
\label{prop:shorten_path}
If $G$ is a graph with
a path $X$ of length 4 whose internal vertices are of degree two and whose end vertices are distinct,
then $\match(G) \simeq \Sigma \match(G/ X)$, where $G / X$ is the contraction of $X$ to a single edge with endpoints given by the endpoints of $X$.
\end{proposition}

The resulting contraction may have parallel edges.

\begin{proposition}[\cite{BJV}, Proposition 11]
\label{double_edge}
Let $G$ be a graph and $e$ an arbitrary edge in $G$. Consider a graph $G \cup \{x\}$ obtained by adding an edge $x$ parallel to $e$ ($x$ and $e$ have same endpoints). Then:
\[\match(G \cup \{x\}) \simeq \match(G) \vee \Sigma \match(G \setminus \operatorname{EN}_{G}[e]).\]
\end{proposition}

We review some notions from graph theory used in the proof.


A {\em cutpoint} of a graph is a vertex whose removal increases the number of components.
A {\em block} of a graph $G$ is a maximal connected subgraph $H$ of $G$
such that $H$ has no cutpoints.
Thus a block with three or more vertices is biconnected.
\begin{definition}
The {\em block-cutpoint graph} of a graph $G$ is the bipartite graph $H_G$
with the following vertex and edge sets.
\begin{itemize}
\item $V(H_G) = \mathcal{B}\cup C$, where $\mathcal{B}$ is the set of blocks of $G$ and $C$ is the set of cut-vertices of $G$;
\item $E(H_G)$ is the set of edges $\{B,c\}$, where $c$ is a cutpoint of $G$
      contained in block $B$ of $G$.
\end{itemize}
\end{definition}

If $G$ is connected, the block-cutpoint graph of $G$ is a tree
whose leaves are blocks of $G$, which we identify as extreme blocks. If $G$ has more than one block, then $H_G$ has at least two leaves, so $G$ has at least two blocks that each contain exactly one cutpoint of $G$ \cite{West}*{p. 156}.

In order to prove our main result, we need the characterization of weak duals of outerplanar graphs.  The \emph{weak dual} of a plane graph is the graph whose vertices correspond to the bounded regions, and two regions are connected if and only if they share one or more edges.
For outerplanar $G$, we assume $G$ is embedded with all vertices on the unbounded face.
It is well-known that the weak dual of an outerplanar graph is a forest, while the weak dual of a biconnected outerplanar graph is a tree \cite{fleischner}.
In an embedded outerplanar graph, the cycles bounding interior regions are called ``basic cycles.'' We use the following result; see Fleischner, Geller, and Harary \cite{fleischner}, Theorem~1, and Sysło \cite{syslo}, Theorem 4.

\begin{theorem}\cite{fleischner,syslo}\label{syslo}
 In a biconnected, outerplanar graph
every pair of basic
cycles has at most one edge in common.
    \end{theorem}

We now prove our main theorem, Theorem~\ref{main-th}.

\begin{proof}
First note that if $G$ is disconnected, so that $G$ is the disjoint union of
two graphs $G_1$ and $G_2$, then $\match(G)=\match(G_1)*\match(G_2)$, and the property of being contractible or homotopy equivalent to a wedge of spheres is preserved by join.
Thus we can assume $G$ is connected.
Also, clearly, any subgraph of an outerplanar graph is outerplanar.

The proof is by induction on the number of edges of the outerplanar graph $G$ with at least one edge.
If the connected graph $G$ has at most 3 edges, then $G$ is  $K_3$ 
(the complete graph on 3 vertices) or a tree, and $\match(G)$ is   
$K_1\sqcup K_2$ or is 0-dimensional, so $\match(G)$ is contractible or a wedge of 0-spheres.
Assume $m>3$ and for every outerplanar graph with fewer than $m$ edges, $\match(G)$ is contractible or homotopy equivalent to a wedge of spheres.  Let $G$ be a connected outerplanar graph with $n$ vertices and $m$ edges.

Suppose $G$ has a vertex $v$ of degree 1.
Let $e$ be the edge containing $v$, and let $h$ be an edge adjacent to $e$.
Then $\operatorname{EN}[e]\subseteq \operatorname{EN}[h]$, so by Proposition~\ref{prop:closed_neigh},
$ \match(G) \simeq \match(G  \setminus\{h \}) \vee \Sigma\match(G  \setminus \operatorname{EN}[h])$.
The graph $G\setminus\{h\}$ is either the disjoint union of $K_2$ (the edge $e$) with an outerplanar graph with two fewer edges than $G$, or is a connected
outerplanar graph with one fewer edge, so its matching complex is,
by induction, contractible or homotopy equivalent to a wedge of spheres.
Also, the graph $G\setminus \operatorname{EN}[h]$ is outerplanar and has fewer edges than $G$, so its matching complex is contractible or homotopy equivalent to a wedge of spheres.  Thus the same holds for $G$.

Now assume $G$ has no vertex of degree 1.
Let $B$ be a block of $G$ that is a leaf of the block-cutpoint graph $H_G$
(an ``extreme block'').
(It might be that $G$ has no cutpoints, in which case $B=G$.)
The subgraph $B$ is itself a biconnected outerplanar graph.
Assume for now that $B$ is not just a single cycle.
The weak dual of $B$ is a tree with at least two vertices, and hence has
at least two leaves.
By Theorem~\ref{syslo}, the vertices on an edge of this tree correspond to a pair of basic cycles in $B$ that share only one edge of $B$.
So the cycles of $B$ corresponding to leaves of the weak dual have the
property that all their vertices, except two adjacent vertices, have degree
2 in $B$.  Let $C_1$ and $C_2$ be two such cycles.
Since $B$ corresponds to a leaf of the block-cutpoint graph,
$B$ contains at most one cutpoint of $G$.  So for at least one of the two
cycles, $C_1$ and $C_2$, the vertices of degree 2 in $B$ are of degree 2
in the whole graph $G$.
We focus on this cycle,
and the path $X$ obtained by removing the edge it shares with the unique adjacent cycle.
We consider different cases
depending on the length of this path.
(See Figure~\ref{figure1}.)

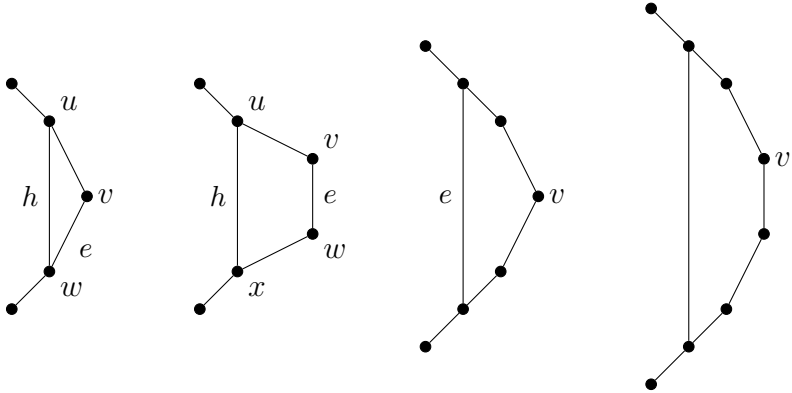
\begin{figure}
    \centering

\begin{tikzpicture}[scale=0.50]
\node[vertex] (a) at (0,2) {};
\node[vertex] (b) at (1,3) {};
\node[vertex] (c) at (2,5) {};
\node[vertex] (d) at (1,7) {};
\node[vertex] (e) at (0,8) {};

\draw (2,5) node[right] {$v$};
\draw (1,3)
node[below right] {$w$};
\draw (1,7) node[above right] {$u$};
\draw (a)--(b);
\draw (c)--(d)--(e);
\draw (b)--(d) node[draw=none,fill=none,font=\small,midway,left] {$h$};
\draw (b)--(c) node[draw=none,fill=none,font=\small,midway,below right] {$e$};

\node[vertex] (f) at (5,2) {};
\node[vertex] (g) at (6,3) {};
\node[vertex] (h) at (8,4) {};
\node[vertex] (i) at (8,6) {};
\node[vertex] (j) at (6,7) {};
\node[vertex] (k) at (5,8) {};

\draw (6,3) node[below right] {$x$};
\draw (8,4) node[below right] {$w$};
\draw (8,6) node[above right] {$v$};
\draw (6,7) node[above right] {$u$};

\draw (f)--(g)--(h);
\draw (h)--(i) node[draw=none,fill=none,font=\small,midway,right] {$e$};
\draw (i)--(j)--(k);
\draw (g)--(j) node[draw=none,fill=none,font=\small,midway,left] {$h$};

\node[vertex] (l) at (11,1) {};
\node[vertex] (m) at (12,2) {};
\node[vertex] (n) at (13,3) {};
\node[vertex] (o) at (14,5) {};
\node[vertex] (p) at (13,7) {};
\node[vertex] (q) at (12,8) {};
\node[vertex] (r) at (11,9) {};

\draw (14,5) node[right] {$v$};

\draw (l)--(m)--(n)--(o)--(p)--(q)--(r);
\draw (m)--(q) node[draw=none,fill=none,font=\small,midway,left] {$e$};

\node[vertex] (s) at (17,0) {};
\node[vertex] (t) at (18,1) {};
\node[vertex] (u) at (19,2) {};
\node[vertex] (v) at (20,4) {};
\node[vertex] (w) at (20,6) {};
\node[vertex] (x) at (19,8) {};
\node[vertex] (y) at (18,9) {};
\node[vertex] (z) at (17,10) {};

\draw (20,6) node[right] {$v$};

\draw (s)--(t)--(u)--(v)--(w)--(x)--(y)--(z);
\draw (t)--(y);
\end{tikzpicture}
\caption{Cases used in proof of Theorem 2}
\label{figure1}
\end{figure}

Case 1.   Suppose the path $X$ has length 2, that is, the two neighbors $u$ and $w$ of $v$ form an edge of $G$.
Let $e=vw$ and $h=uw$.  Then the closed neighborhood $\operatorname{EN}[e]$ is contained in
the closed neighborhood $\operatorname{EN}[h]$, so
$ \match(G) \simeq {\match(G  \setminus\{h \})} \vee {\Sigma\match(G  \setminus \operatorname{EN}[h])}$, and, as in the case of a leaf, by induction, $G$ has matching complex contractible or homotopy equivalent to a wedge of spheres.

Case 2.  Suppose the path $X$ has length 3.  Then $G$ contains an induced
4-cycle $u,v,w,x$, where $w$ also has degree 2 in $G$.
Let $e=vw$, $h=ux$.  The open neighborhood of $e$ is
$\operatorname{EN}(e)= \{uv, wx\}$; this set is contained in the open neighborhood of $\operatorname{EN}(h)$.
So by Proposition~\ref{prop:open_neighborhood},
$\match(G)\simeq\match(G\setminus\{h\})$.  By induction, $\match(G)$ is
contractible or homotopy equivalent to a wedge of spheres.

Case 3.  Suppose the path $X$ has length 4.
Let $e$ be the edge of $G$
connecting the two endpoints of $X$.
Then by Proposition~\ref{prop:shorten_path},
$\match(G)\simeq\Sigma\match(G/X)$, where $G/X$ is the contraction of $X$
to a single edge $x$.  
Let $G\setminus X$ be the graph obtained by removing all edges and internal 
vertices of the path $X$ from $G$, and note that 
$G/X= (G\setminus X)\cup \{x\}$.
Apply Proposition~\ref{double_edge} to the graph $G\setminus X$ with edge $e$.
So 
$$\match(G/X) = {\match((G\setminus X)\cup\{x\})}={\match(G\setminus X)\vee 
\Sigma\match((G\setminus X)\setminus \operatorname{EN}_{G\setminus X}[e])}.$$ 
Each of $G\setminus X$ and $(G\setminus X)\setminus \operatorname{EN}_{G\setminus X}[e]$ is outerplanar with
fewer edges than $G$, so by induction, $\match(G/X)$ and, hence, $\match(G)$ is contractible or homotopy equivalent to a wedge of spheres.

Case 4.  Suppose the path $X$ has length 5 or more.  Let $Y$ be a subpath
of length 4.  By
Proposition~\ref{prop:shorten_path},
$\match(G)\simeq\Sigma\match(G/Y)$, where $G/Y$ is the contraction of $Y$
to a single edge $x$. In this case $G/Y$ is an outerplanar graph with fewer edges than
$G$, so by induction,
$\match(G)$ is contractible or homotopy equivalent to a wedge of spheres.

In these arguments we have assumed that $B$ is not a single cycle.  If $B$ is a single cycle of length $n$, then similar arguments, using Proposition~\ref{prop:closed_neigh} for a cycle of length 3, Proposition~\ref{prop:open_neighborhood} for a cycle of length 4, and Proposition~\ref{prop:shorten_path} for a cycle of length 5 or more, allow us to reduce the matching
complex to matching complexes of outerplanar graphs with fewer edges.
So by induction $\match(G)$ is contractible or homotopy equivalent to a wedge of spheres.
\end{proof}

\section{Conclusion}

We were inspired to study the homotopy type of matching complexes of outerplanar graphs based on recent results showing some outerplanar graphs have matching complexes that are contractible or homotopy equivalent to wedges of spheres.  These include
paths and cycles \cite{kozlov-directed}, forests \cite{MT_forests, JelicEtAl}, and polygonal line tilings \cite{takahiro,JelicEtAl}.

The next step is to see whether this can be extended to planar graphs.  One class of planar (nonouterplanar) graphs whose matching complexes are known to be  contractible or homotopy equivalent to a wedge of spheres is the categorical product of paths \cite{gupta2024matching}.
Other examples are wheel graphs \cite{Vega}, planar complete bipartite graphs (namely, $K_{m,2}$) and the graphs formed by joining copies of $K_4$ at a common vertex.  We now have a new tool for proving matching complexes are contractible or homotopy equivalent to wedges of spheres: use Propositions~\ref{prop:math_ind}, \ref{prop:closed_neigh}, \ref{prop:open_neighborhood}, and~\ref{prop:shorten_path} to reduce the graph to an outerplanar graph.  

We do not know of an example of a planar graph whose matching complex is neither contractible nor homotopy equivalent to a wedge of spheres.  Hence, we propose the following. 

\begin{conjecture} The matching complex of every planar graph is contractible or homotopy equivalent to a wedge of spheres.
\end{conjecture}

Besides homotopy type, the homology, and in particular, the existence of torsion in homology, of matching complexes is of interest.  Torsion does appear in the matching complexes of complete graphs and complete bipartite graphs (Wachs \cite{wachs}, Jonsson \cite{jonsson-book}). We end with a final question:

 Which graph properties give rise to torsion in matching complexes?

\end{document}